\documentclass[oneside,10pt]{article}          
\usepackage[b5paper]{geometry}	    
\usepackage{amsfonts,amsmath,latexsym,amssymb} 
\usepackage{amsthm}                
\usepackage{mathrsfs,upref}         
\usepackage{mathptmx}
\usepackage{leftidx}
\usepackage{mathtools}	
\usepackage{enumitem}	    
\usepackage{dea}	            
\newtheorem{thm}{Theorem}[section]
\newtheorem{cor}{Corollary}[section]
\newtheorem{lem}{Lemma}[section]
\newtheorem{defn}{Definition}[section]
\newtheorem{rem}{Remark}[section]
\newtheorem{proposition}{Proposition}[section]

\newcounter{romnum}

\begin{document}

\title[Lyapunov-type Inequalities]
{Lyapunov-Type Inequalities for Third Order Nonlinear Equations}

\author{Brian Behrens and Sougata Dhar}

\address{Brian Behrens, Department of Mathematics, University of Connecticut, Storrs, CT 06269, USA.\\
\email{brian.behrens@uconn.edu}}

\address{Sougata Dhar, Department of Mathematics, University of Connecticut, Storrs, CT 06269, USA.\\
\email{brian.behrens@uconn.edu}}

\CorrespondingAuthor{Sougata Dhar}

\dedicated{Dedicated to Prof. Paul Eloe}                    

\date{09.02.2022}                               



\begin{abstract} 
We derive Lyapunov-type inequalities for general third order nonlinear equations involving multiple $\psi$-Laplacian operators of the form 
\begin{equation*}
    (\psi_{2}((\psi_{1}(u'))'))' + q(x)f(u) = 0,
\end{equation*} 
where $\psi_{2}$ and $\psi_{1}$ are odd, increasing functions, $\psi_{2}$ is super-multiplicative, $\psi_{1}$ is sub-multiplicative, and $\frac{1}{\psi_{1}}$ is convex, and $f$ is a continuous function which satisfies a sign condition. Our results utilize $q_{+}$ and $q_{-}$, as opposed to $|q|$ which appears in most results in the literature. Additionally, these new inequalities generalize previously obtained results, and the proofs utilize a different technique than most other works in the literature. Furthermore, using the obtained inequalities, we obtain a constraint on the location of the maximum of a solution, properties of oscillatory solutions, and an upper bound for the number of zeroes.  
\end{abstract}

\maketitle



\setcounter{equation}{0}
\section{Introduction}

In this paper, we establish Lyapunov-type inequalities for a generalized third order nonlinear differential equation in the form
\begin{equation}\label{1.1}
    (\psi_{2}((\psi_{1}(u'))'))' + q(x)f(u) = 0,
\end{equation} 
with one of the following two boundary conditions:
\begin{align}
    u(a) =  u(b) = 0, u(x) \neq 0 \text{ for } x \in (a,b), (\psi_{1}(u'))'(\xi) = 0 \text{ for some } \xi \in [a,b], \label{BC1}
\end{align} 
and
\begin{align}
    u(a) =  u(b) = u(c) = 0, u(x) \neq 0 \text{ for } x \in (a,b)\cup (b,c). \label{BC2}
\end{align}

\newpage

Throughout this paper, unless mentioned otherwise, we will assume the following:
\begin{enumerate}
    \item[(H1)] $\psi_{1}$ and $\psi_{2}$ are odd and increasing functions.
    \item[(H2)] $\psi_{1}(s)$ is sub-multiplicative on $[0,\infty)$ and $\frac{1}{\psi_{1}(s)}$ is a convex function for $s>0$.
    \item[(H3)] $\psi_{2}(s)$ is super-multiplicative on $[0,\infty)$
    \item[(H4)] $q\in C([a,b],\mathbb{R})$ and can change sign in $[a,b]$. 
    \item[(H5)] $f:\mathbb{R}\to\mathbb{R}$ is odd continuous function and satisfies $sf(s) > 0$ for $s \neq 0$. 
\end{enumerate}

To establish the Lyapunov-type inequalities for a general nonlinear equation \eqref{1.1}, it is fundamental that $\psi_{1}$ is sub-multiplicative and $\psi_{2}$ is super-multiplicative. Below we provide the definition of a 
sub and super-multiplicative functions, respectively.

\begin{defn}
 A function $v:A \to [0,\infty]$ is sub-multiplicative if $v(xy) \leq v(x)v(y)$ for all $x,y \in A$. A function $w:A \to [0,\infty]$ is super-multiplicative if $w(xy)\ge w(x)w(y)$ for all $x,y \in A$.
\end{defn}

In \cite{N4, R5}, the basic properties and some examples of sub-multiplicative functions are given. These class of functions
frequently appear in the semi-group theory \cite{R6}, interpolation theory \cite{N1, N2} play a significant role in the theory of 
operators on Orlicz spaces, see \cite{N3} for details.
The power function $v(x)=x^p\, (x>0)$, where $p\in\mathbb{R}$ is a common example of the sub-multiplicative class of 
functions. This gives rise to the classical $p-$Laplacian $\phi_p(x)=|x|^{p-1}x$, $p>1$ as an important example 
of a sub-multiplicative function.
Note that Eq.~\eqref{1.1} becomes an equation involving the classical $p$-Laplacian operator when $\psi_{1}$ and $\psi_{2}$
are replaced by $\phi_p$'s. A similar discussion can be made about super-multiplicative functions. In fact, the motivation behind this paper is that the sub-multiplicative and super-multiplicative class of functions
allow us to consider a much greater class of operators and nonlinearities which in turn provides the most general result
in the literature. In particular, our proofs do not require $\psi_{1}$, $\psi_{2}$ and $f$ to be power functions but covers those
as special cases. 

Now we briefly visit the historical developments of the Lyapunov-type inequalities.
In \cite{R1}, A.M. Lyapunov considered the following second order linear boundary value problem (BVP)
\begin{equation}\label{1.2}
    u'' + q(x)u = 0,\quad  u(a) = u(b) = 0,
\end{equation} 
and obtained the following result.
\begin{thm}\label{t1}
Assume \eqref{1.2} has a nontrivial solution $u(x)$ and $u(x)\neq 0$ for $x\in(a,b)$. Then
\begin{equation}\label{1.3}
    \int_{a}^{b}|q(x)|dx > \frac{4}{b-a}.
\end{equation}
\end{thm} 
The inequality \eqref{1.3} is commonly referred to as the Lyapunov inequality. It has wide reaching applications in many areas 
of differential equations including eigenvalue problems, boundary value problems, and oscillation to name a few. 
As a result, it has been improved and extended in numerous directions. 
The first improvement on Lyapunov's result was done by Wintner \cite{R2}, by replacing $|q(x)|$ in \eqref{1.3} with $q_{+}(x) := \max_{x\in [a,b]}\{q(x),0\}$. In particular, under the same assumptions as in Theorem \ref{t1}, he obtained 
\begin{equation}\label{1.4}
    \int_{a}^{b}q_{+}(x)dx > \frac{4}{b-a}.
\end{equation}
Yang \cite{R3} extended \eqref{1.4} in the direction of a second order quasilinear BVP
\begin{equation*}
         (\phi_{\alpha}(u'))' + q(x)\phi_{\alpha}(u) = 0,\quad     u(a) = u(b) = 0,
\end{equation*} 
where $\phi_{\alpha}(u) = |u|^{\alpha - 1}u, \,\alpha > 0$, and obtained 
\begin{equation}\label{1.5}
    \int_{a}^{b}q_{+}(x)dx > \frac{2^{\alpha +1}}{(b-a)^{\alpha}}.
\end{equation}
It is clear that \eqref{1.5} reduces to \eqref{1.4} when $\alpha=1$. 
Further improvements and generalizations of Lyapunov-type inequalities have been made by a number of researchers and we 
refer the interested reader to \cite{S2, T5, T7, S1, T3, S3, T6, T1, R4}.
 Among those, a particularly interesting result 
was formulated by Sanchez and Vergara \cite{R4}, where the authors considered a general nonlinear second order equation
\begin{equation}\label{1.6}
         (\psi(u'))' + q(x)f(u) = 0.
\end{equation} 
Here $\psi(u)$ is a sub-multiplicative operator. Under the same assumptions as in Theorem \ref{t1}, the authors obtained
\begin{equation}\label{1.7}
    \int_{a}^{b}q_{+}(x)\frac{f(u)}{\psi(u)}\,dx > \frac{2}{\psi(\frac{b-a}{2})}.
\end{equation}
Inequality \eqref{1.7} is the most general for the second order case as it includes \eqref{1.5} and \eqref{1.4} for some specific
choices of the operator $\psi$. Although, the Lyapunov-type inequalities are well developed
for the even order equations, there were no results existed for the odd order equations until late $1990$. It was 
Parhi and Panigrahi \cite{R7}, who first developed the Lyapunov-type inequalities in 1999 for the third order linear equation \begin{equation}\label{1.8}
    u''' + q(x)u = 0,
\end{equation} 
with each of the following two boundary conditions.
\begin{equation}
     u(a) =  u(b) = 0 \text{ and } u''(\xi) = 0 \text{ for some } \xi\in[a,b], \label{BCPP1}
\end{equation}
and
\begin{equation}
     u(a) =  u(b) = u(c) = 0. \label{BCPP2}
\end{equation}
We summarize their results in the following theorem. 
\begin{thm} \label{t2}
(a) Assume Eq.~\eqref{1.8} has a nontrivial solution $u(x)$ satisfying \eqref{BCPP1}. Then
\begin{equation}\label{1.9}
    \int_{a}^{b}|q(x)|dx >\frac{4}{(b-a)^2}.
\end{equation}

(b) Assume Eq.~\eqref{1.8} has a nontrivial solution $u(x)$ satisfying \eqref{BCPP2}. Then 
\begin{equation}\label{1.10}
    \int_{a}^{c}|q(x)|dx >\frac{4}{(c-a)^2}.   
\end{equation} 
\end{thm} 
The work of Parhi and Panigrahi inspired researchers to establish Lyapunov-type inequalities for odd order cases, 
as seen in \cite{T9, R8, R9, R10, R11, T2, T8, T10, T11}. 
The third order quasilinear equation was studied by Dhar and Kong in $2014$. In \cite{R8}, the authors considered
\begin{equation}\label{1.11}
    (\phi_{\alpha_{2}}((\phi_{\alpha_{1}}(u'))'))' + q(x)\phi_{\alpha_{1}\alpha_{2}}(u) = 0,
\end{equation} 
where $\phi_{p}(x) = |x|^{p-1}x$, $p> 0$ and $p=\alpha_{1},\alpha_{2}$
with each of the following two boundary conditions.
\begin{align}
    u(a) =  u(b) = 0, u(x) \neq 0 \text{ for } x \in (a,b), (\phi_{\alpha_{1}}(u'))'(\xi) = 0 \text{ for some } \xi\in[a,b] ,\label{BCDK1}
\end{align} 
and
\begin{align}
    u(a) =  u(b) = u(c) = 0, u(x) \neq 0 \text{ for } x \in (a,b)\cup (b,c). \label{BCDK2}
\end{align}
We summarize their results in the following theorem. 
\begin{thm}
(a) Assume Eq.~\eqref{1.11} has a nontrivial solution $u(x)$ satisfying \eqref{BCDK1}. Then 
\begin{equation}
    \int_{a}^{\xi}q_{-}(x)\,dx + \int_{\xi}^{b}q_{+}(x)\,dx > \left(\frac{2}{b-a}\right)^{(\alpha_{1}+1)\alpha_{2}}.
\end{equation} 

(b) Assume Eq.~\eqref{1.11} has a nontrivial solution $u(x)$ satisfying \eqref{BCDK2}. Then 
\begin{equation}
    \max_{\xi\in[a,c]}\left\{\int_{a}^{\xi}q_{-}(x)\,dx + \int_{\xi}^{b}q_{+}(x)\,dx\right\} > \left(\frac{2}{c-a}\right)^{(\alpha_{1}+1)\alpha_{2}}.
\end{equation}
\end{thm} 
Their result is also an improvement of Theorem \ref{t2} even for the linear case, using the positive $(q_+(x))$ and negative 
$(q_{-}(x) := \max_{x\in [a,b]}\{-q(x),0\})$ parts of $q(x)$, as opposed to $|q(x)|$. 

In this article, we obtain the Lyapunov-type inequalities for the third order nonlinear equations in the form 
\eqref{1.1}. Our work is an extension of Sanchez and Vergara's work \cite{R4} from the second order case to the third order 
case. Moreover, our results cover the results obtained by Dhar and Kong \cite{R8} as a special case and provides the
most general result for the third order nonlinear equations. In addition to that, we employ new techniques to prove 
our main results which does not require the use of the Cauchy-Schwartz inequality as is seen in the majority of the works in the literature. As applications of our obtained results, we discuss the upper bound on the number of zeros of a 
nontrivial solution, criterion for oscillatory solutions, and comment on the location of the maximum of a nontrivial solution.

\setcounter{equation}{0}
\section{Main Results}
We first recall some well known results from the literature which will be useful to prove our main theorems and subsequent results. 
\begin{lem}[Jensen's Inequality]\label{lemma2.1} 
Let $g$ be a convex function and $t \in [0,1]$. Then
\begin{equation}\label{2.2}
g(tx_{1} + (1-t)x_{2}) \le tg(x_{1}) + (1-t)g(x_{2}).
\end{equation} 
Furthermore, 
\begin{equation}\label{2.3}
 g\left(\frac{1}{N}\sum_{k=1}^{N}x_{k}\right) \leq \frac{1}{N}\sum_{k=1}^{N}g(x_{k}), \quad N \in \mathbb{N}.
\end{equation}
\end{lem}

We are now prepared to present our first result on Lyapunov-type Inequalities.

\begin{thm}
\label{thm2.1} Assume $(H1)-(H5)$ holds and Eq.~\eqref{1.1} has a nontrivial solution $u(x)$ satisfying \eqref{BC1}.
Then
\begin{equation}
\int_{a}^{\xi}q_{-}(x)\Phi(u)\,dx + \int_{\xi}^{b}q_{+}(x)\Phi(u)\,dx > \psi_{2}\left(\frac{\tfrac{2}{b-a}}{\psi_{1}\left(\tfrac{b-a}{2}\right)}\right),\label{2.4}
\end{equation}
where $\Phi(u) = \frac{f(u)}{\psi_{2}(\psi_{1}(u))}$.
\end{thm}
\begin{proof}
Without loss of generality, assume $u(x) > 0$ on $(a,b)$. Since $u(x)$ satisfies \eqref{BC1}, there exists $c \in (a,b)$ such that $u(c) = \max_{x\in[a,b]} u(x)$. By the Mean Value Theorem, there exist $\tau_{1} \in (a,c)$ and $\tau_{2} \in (c,b)$ such that \begin{align}\label{2.5}
u'(\tau_{1}) = \frac{u(c) - u(a)}{c-a} = \frac{u(c)}{c-a} \text{ and }
u'(\tau_{2}) = \frac{u(b) - u(c)}{b-c} = \frac{-u(c)}{b-c}.
\end{align} 
Recall that $\frac{1}{\psi_{1}(x)}$ is convex. Hence using \eqref{2.2} with $t = \tfrac{1}{2}$, $x_{1} = c-a$, and $x_{2} = b-c$, we have \begin{equation*}
    \frac{2}{\psi_{1}\left(\tfrac{b-a}{2}\right)} \leq \frac{1}{\psi_{1}(c-a)} + \frac{1}{\psi_{1}(b-c)} = \frac{1}{\psi_{1}(u(c))}\left[\frac{\psi_{1}(u(c))}{\psi_{1}(c-a)} + \frac{\psi_{1}(u(c))}{\psi_{1}(b-c)}\right].
\end{equation*} 
It follows that 
\begin{eqnarray}
    \frac{2\psi_{1}(u(c))}{\psi_{1}\left(\tfrac{b-a}{2}\right)} &\leq& \frac{\psi_{1}\left((c-a)\frac{u(c)}{c-a}\right)}{\psi_{1}(c-a)} + \frac{\psi_{1}\left((b-c)\frac{u(c)}{b-c}\right)}{\psi_{1}(b-c)}. \nonumber 
    \end{eqnarray} 
Since $\psi_{1}$ is sub-multiplicative, we obtain 
\begin{eqnarray}
    \frac{2\psi_{1}(u(c))}{\psi_{1}\left(\tfrac{b-a}{2}\right)} &\leq& \psi_{1}\left(\frac{u(c)}{c-a}\right) + \psi_{1}\left(\frac{u(c)}{b-c}\right). \nonumber 
    \end{eqnarray} 
    Using (\ref{2.5}), we have
    \begin{eqnarray}
    \frac{2\psi_{1}(u(c))}{\psi_{1}\left(\tfrac{b-a}{2}\right)} \leq \psi_{1}(u'(\tau_{1})) - \psi_{1}(u'(\tau_{2})) = \int_{\tau_{1}}^{\tau_{2}} \label{2.6} -(\psi_{1}(u'))'\,dx.
\end{eqnarray}
Recall that $(\psi_{1}(u'))'(\xi) = 0$ and $\psi_{2}$ is odd, we have $\psi_{2}((\psi_{1}(u'))')(\xi) = 0$. Integrating \eqref{1.1} from $\xi$ to $t$ and using this fact, we have
\begin{equation*}
    \psi_{2}((\psi_{1}(u'))') = -\int_{\xi}^{x}q(s)f(u)\,ds.
\end{equation*} 
Since $\psi_{2}$ is odd, it follows that  
\begin{equation*}
     (\psi_{1}(u'))' = \psi_{2}^{-1}\left(-\int_{\xi}^{x}q(s)f(u)\,ds\right) = -\psi_{2}^{-1}\left(\int_{\xi}^{x}q(s)f(u)\,ds\right).
\end{equation*} 
Replacing $(\psi_{1}(u'))'$ in (\ref{2.6}), we obtain 
\begin{equation*}
    \frac{2\psi_{1}(u(c))}{\psi_{1}\left(\tfrac{b-a}{2}\right)} \leq \int_{a}^{b}\psi_{2}^{-1}\left(\int_{\xi}^{x}q(s)f(u)\,ds\right)\,dx.
\end{equation*} 
Hence
\begin{eqnarray*}
    \frac{2\psi_{1}(u(c))}{\psi_{1}\left(\tfrac{b-a}{2}\right)} 
    &\leq& \int_{a}^{\xi}\psi_{2}^{-1}\left(\int_{\xi}^{x}q(s)f(u)\,ds\right) \,dx + \int_{\xi}^{b} \psi_{2}^{-1}\left(\int_{\xi}^{x}q(s)f(u)\,ds\right) \,dx \\
    &=& \int_{a}^{\xi}\psi_{2}^{-1}\left(\int_{x}^{\xi}-q(s)f(u)\,ds\right) \,dx + \int_{\xi}^{b} \psi_{2}^{-1}\left(\int_{\xi}^{x}q(s)f(u)\,ds\right) \,dx.
    \end{eqnarray*} 
Since $-q(x)\le q_{-}(x)$ and $ q(x)\le q_{+}(x)$, we have 
\begin{eqnarray*}
    \frac{2\psi_{1}(u(c))}{\psi_{1}\left(\tfrac{b-a}{2}\right)} 
    &\le& \int_{a}^{\xi} \psi_{2}^{-1}\left(\int_{a}^{\xi}q_{-}(s)f(u)\,ds\right) \,dx  + \int_{\xi}^{b} \psi_{2}^{-1}\left(\int_{\xi}^{b}q_{+}(s)f(u)\,ds\right) \,dx \\
    &=& (\xi-a)\psi_{2}^{-1}\left(\int_{a}^{\xi}q_{-}(s)f(u)\,ds\right)+(b-\xi)\psi_{2}^{-1}\left(\int_{\xi}^{b}q_{+}(s)f(u)\,ds\right).
    \end{eqnarray*} 
Dividing both sides by $b - a$, we have  
\begin{equation}\label{2.7}
    \frac{\frac{2}{b-a}\psi_{1}(u(c))}{\psi_{1}\left(\tfrac{b-a}{2}\right)} \leq \frac{\xi-a}{b-a}\psi_{2}^{-1}\left(\int_{a}^{\xi}q_{-}(s)f(u)\,ds\right)+\frac{b-\xi}{b-a}\psi_{2}^{-1}\left(\int_{\xi}^{b}q_{+}(s)f(u)\,ds\right).
\end{equation} 

We discuss two cases based on the nature of $\psi_{2}$. \\ 
Case 1: Assume $\psi_{2}$ is convex. Applying $\psi_{2}$ to (\ref{2.7}), we have \begin{eqnarray}
    \psi_{2}\left(\frac{\tfrac{2}{b-a}\psi_{1}(u(c))}{\psi_{1}\left(\tfrac{b-a}{2}\right)}\right) \leq \psi_{2}\biggl(\frac{\xi-a}{b-a}\psi_{2}^{-1}\left(\int_{a}^{\xi}q_{-}(x)f(u)\,dx\right) + \frac{b-\xi}{b-a}\psi_{2}^{-1}\left(\int_{\xi}^{b}q_{+}(x)f(u)\,dx\right)\biggl). \nonumber \end{eqnarray} 
Using \eqref{2.2} with $t=\frac{\xi-a}{b-a}$, we have 
\begin{eqnarray}
    \psi_{2}\left(\frac{\tfrac{2}{b-a}\psi_{1}(u(c))}{\psi_{1}\left(\tfrac{b-a}{2}\right)}\right) &\leq& \frac{\xi-a}{b-a}\int_{a}^{\xi}q_{-}(x)f(u)\,dx + \frac{b-\xi}{b-a}\int_{\xi}^{b}q_{+}(x)f(u)\,dx, \nonumber \\
    &\leq& \int_{a}^{\xi}q_{-}(x)f(u)\,dx + \int_{\xi}^{b}q_{+}(x)f(u)\,dx. \nonumber
\end{eqnarray} \\
Case 2: Assume $\psi_{2}$ is concave. Then for $x_{1},x_{2} \in [0, \infty)$, we have \begin{equation}
    \psi_{2}(x_{1} + x_{2}) \leq \psi_{2}(x_{1}) + \psi_{2}(x_{2}). \label{concavity}
\end{equation} 
Recall that $\frac{\xi-a}{b-a}, \frac{b-\xi}{b-a} \leq 1$, so from (\ref{2.7}) we have \begin{equation*}
    \frac{\frac{2}{b-a}\psi_{1}(u(c))}{\psi_{1}\left(\tfrac{b-a}{2}\right)} \leq \psi_{2}^{-1}\left(\int_{a}^{\xi}q_{-}(s)f(u)\,ds\right)+\psi_{2}^{-1}\left(\int_{\xi}^{b}q_{+}(s)f(u)\,ds\right).
\end{equation*} Applying $\psi_{2}$ and using (\ref{concavity}), we have \begin{equation*}
    \psi_{2}\left(\frac{\tfrac{2}{b-a}\psi_{1}(u(c))}{\psi_{1}\left(\tfrac{b-a}{2}\right)}\right) \leq \int_{a}^{\xi}q_{-}(x)f(u)\,dx + \int_{\xi}^{b}q_{+}(x)f(u)\,dx.
\end{equation*} In each case, we have \begin{equation*}
    \psi_{2}\left(\frac{\tfrac{2}{b-a}\psi_{1}(u(c))}{\psi_{1}\left(\tfrac{b-a}{2}\right)}\right) \leq \int_{a}^{\xi}q_{-}(x)f(u)\,dx + \int_{\xi}^{b}q_{+}(x)f(u)\,dx.
\end{equation*} Using the fact that $\psi_{2}$ is super-multiplicative, we have 
\begin{equation*}
    \psi_{2}(\psi_{1}(u(c)))\psi_{2}\left(\frac{\tfrac{2}{b-a}}{\psi_{1}\left(\tfrac{b-a}{2}\right)}\right) \leq \int_{a}^{\xi}q_{-}(x)f(u)\,dx + \int_{\xi}^{b}q_{+}(x)f(u)\,dx.
\end{equation*} Therefore because $\psi_{1}$ and $\psi_{2}$ are increasing, we obtain \begin{equation*}
     \psi_{2}\left(\frac{\tfrac{2}{b-a}}{\psi_{1}\left(\tfrac{b-a}{2}\right)}\right) < \int_{a}^{\xi}q_{-}(x)\Phi(u)\,dx + \int_{\xi}^{b}q_{+}(x)\Phi(u)\,dx,
\end{equation*} where $\Phi(u) = \frac{f(u)}{\psi_{2}(\psi_{1}(u))}$. The proof is now complete.
\end{proof}

In the next theorem, we present a Lyapunov-type Inequality for Eq.~\eqref{1.1} with a three point boundary condition
\eqref{BC2}. 

 \begin{thm}\label{t2.2a}
\label{thm2.2} Assume $(H1)-(H5)$ holds and Eq.~\eqref{1.1} has a nontrivial solution $u(x)$ satisfying \eqref{BC2}. Then either 
\begin{equation}\label{2.9}
\max_{\xi\in [a,b]}\left\{\int_{a}^{\xi}q_{-}(x)\Phi(u)\,dx + \int_{\xi}^{b}q_{+}(x)\Phi(u)\,dx\right\} > \psi_{2}\left(\frac{\tfrac{2}{b-a}}{\psi_{1}\left(\tfrac{b-a}{2}\right)}\right),
\end{equation} 
or 
\begin{equation}\label{2.10}
\max_{\xi\in [b,c]}\left\{\int_{b}^{\xi}q_{-}(x)\Phi(u)\,dx + \int_{\xi}^{c}q_{+}(x)\Phi(u)\,dx\right\} > \psi_{2}\left(\frac{\tfrac{2}{c-b}}{\psi_{1}\left(\tfrac{c-b}{2}\right)}\right).
\end{equation} 
As a result, 
\begin{equation}\label{2.11}
\max_{\xi\in [a,c]}\left\{\int_{a}^{\xi}q_{-}(x)\Phi(u)\,dx + \int_{\xi}^{c}q_{+}(x)\Phi(u)\,dx\right\} > \psi_{2}\left(\frac{\tfrac{2}{c-a}}{\psi_{1}\left(\tfrac{c-a}{2}\right)}\right),
\end{equation} where $\Phi(u) = \frac{f(u)}{\psi_{2}(\psi_{1}(u))}$.
\end{thm}

\begin{proof}
By Rolle's Theorem, there exist $x_{1} \in [a,b]$ and $x_{2} \in [b,c]$ such that $u'(x_{1}) = u'(x_{2}) = 0$. Further application of the Mean Value Theorem shows that there exists some $\xi \in (x_{1}, x_{2})$ such that \begin{equation*}
    \left(\psi_{1}\left(u'(\xi)\right)\right)' = \frac{\psi_{1}\left(u'(x_{2})\right) - \psi_{1}\left(u'(x_{2})\right)}{x_{2} - x_{1}} = \frac{\psi_{1}(0) - \psi_{1}(0)}{x_{2} - x_{1}} = 0.
\end{equation*} Since $\psi_{2}$ is odd, we have \begin{equation*}
    \psi_{2}\left(\left(\psi_{1}\left(u'(\xi)\right)\right)'\right) = 0.
\end{equation*} \\
Clearly, either $\xi \in (a,b]$ or $\xi \in [b,c)$. If $\xi \in (a,b]$, then applying Theorem 2.1, we have \begin{equation}
    \psi_{2}\left(\frac{\tfrac{2}{b-a}}{\psi_{1}\left(\tfrac{b-a}{2}\right)}\right) < \int_{a}^{\xi}q_{-}(x)\Phi(u)\,dx + \int_{\xi}^{b}q_{+}(x)\Phi(u)\,dx,
\end{equation} hence \eqref{2.9} holds. Similarly, if $\xi \in [b,c)$ then we have \begin{equation}
    \psi_{2}\left(\frac{\tfrac{2}{c-b}}{\psi_{1}\left(\tfrac{c-b}{2}\right)}\right) < \int_{b}^{\xi}q_{-}(x)\Phi(u)\,dx + \int_{\xi}^{c}q_{+}(x)\Phi(u)\,dx. 
\end{equation} 
Hence \eqref{2.10} holds. It is trivial to see that both \eqref{2.9} and \eqref{2.10} lead to \eqref{2.11}.  
This concludes the proof.
\end{proof}

Since $q_{-}(x), q_{+}(x) \leq |q(x)|$, the following results follow directly from Theorems 2.1 and 2.2. 

\begin{cor}\label{cor2.1} 
(a) Assume $(H1)-(H5)$ holds and Eq.~\eqref{1.1} has a nontrivial solution $u(x)$ satisfying \eqref{BC1}. Then \begin{equation*}
        \int_{a}^{b}|q(x)|\Phi(u)\,dx > \psi_{2}\left(\frac{\tfrac{2}{b-a}}{\psi_{1}\left(\tfrac{b-a}{2}\right)}\right).
    \end{equation*} \\
(b) Assume $(H1)-(H5)$ holds and Eq.~\eqref{1.1} has a nontrivial solution $u(x)$ satisfying \eqref{BC2}. Then \begin{equation*}
   \int_{a}^{b}|q(x)|\Phi(u)\,dx > \psi_{2}\left(\frac{\tfrac{2}{b-a}}{\psi_{1}\left(\tfrac{b-a}{2}\right)}\right), \end{equation*}  and  \begin{equation*} \int_{b}^{c}|q(x)|\Phi(u)\,dx > \psi_{2}\left(\frac{\tfrac{2}{c-b}}{\psi_{1}\left(\tfrac{c-b}{2}\right)}\right).
    \end{equation*}
As a result, \begin{equation*}
    \int_{a}^{c}|q(x)|\Phi(u)\,dx > \psi_{2}\left(\frac{\tfrac{2}{c-a}}{\psi_{1}\left(\tfrac{c-a}{2}\right)}\right).
\end{equation*}
\end{cor}

It is clear that the results provided in Corollary \ref{cor2.1} is simpler that those given in Theorems 2.1 and 2.2. 
However, these are much weaker as well. For example, if $\xi=a$ then Theorem 2.1 lead to
$$
\int_{a}^{b}q_{-}(x)\Phi(u)\,dx >\psi_{2}\left(\frac{\tfrac{2}{b-a}}{\psi_{1}\left(\tfrac{b-a}{2}\right)}\right),
$$
and if $\xi=b$ then Theorem 2.1 lead to
$$
\int_{a}^{b}q_{+}(x)\Phi(u)\,dx > \psi_{2}\left(\frac{\tfrac{2}{b-a}}{\psi_{1}\left(\tfrac{b-a}{2}\right)}\right).
$$
These can not be observed from Corollary \ref{cor2.1}.

\begin{rem}\label{rem2.1}
\rm{
Let $\psi_{1}$, $\psi_{2},$ and $f$ be signed power functions of the form $\psi_{1}(x) = |x|^{\alpha_{1}-1}x$, 
$\psi_{2}(x) = |x|^{\alpha_{2}-1}x$, and $f(x) = |x|^{\alpha_{2}\alpha_{1}-1}x$ for $\alpha_{1}, \alpha_{2} > 0$. 
Hence, \eqref{1.1} reduces to the quasilinear equation \eqref{1.11} which was studied by Dhar and Kong in \cite{R7}. 
In this case, $\Phi(u)=1$ and 
$$
\psi_{2}\left(\frac{\tfrac{2}{b-a}}{\psi_{1}\left(\tfrac{b-a}{2}\right)}\right)=\left(\frac{2}{b-a}\right)^{\alpha_{2}(\alpha_{1} + 1)}.
$$
Thus we have by Theorem 2.1
\begin{equation*}
     \int_{a}^{\xi} q_{-}(x)\,dx + \int_{\xi}^{b} q_{+}(x)\,dx > \left(\frac{2}{b-a}\right)^{\alpha_{2}(\alpha_{1} + 1)},
\end{equation*} 
which is consistent with the results obtained by the authors in \cite{R8}. A same argument holds true for 
$\alpha_{1} = \alpha_{2} = 1$, i.e., the linear case. We leave the details for the interested reader.
}
\end{rem}

Let us now suppose that the nonlinearity of $f$ lies between power functions. In particular, we suppose there exist 
positive constants $p, c_{1}, \text{ and }c_{2}$ such that 
\begin{equation}
    c_{1}|s|^{p} \leq |f(s)| \leq c_{2}|s|^{p} \text{ for all } s \in \mathbb{R}. \label{fineq}
\end{equation}
We are concerned with the location of the point where a solution of Eq.~\eqref{1.1} attains a maximum, 
specifically with the case when $\psi_{1}$ and $\psi_{2}$ are power functions. For the following propositions, we will assume 
 $\psi_{1}(s) = |s|^{\alpha_{1}-1}s$ and $\psi_{2}(s) = |s|^{\alpha_{2}-1}s$, and the inequality in \eqref{fineq} holds.

\begin{proposition}
    Assume $(H1)-(H4)$ holds and Eq.~\eqref{1.1} has a nontrivial solution $u(x)$ satisfying \eqref{BC1}. 
    Let $M = |u(d)| = ||u||_{\infty} = \max_{a\leq x\leq b}|u(x)| > 0$, and suppose $p - \alpha_{1}\alpha_{2} > 0$. 
    Then $d$ cannot be too close to $a$ or $b$. \end{proposition}
    \begin{proof}
    By Corollary 2.1, we have 
    \begin{align*}
    \left(\frac{2}{b-a}\right)^{\alpha_{2}(\alpha_{1} + 1)} &\le \int_{a}^{b}|q(x)|\frac{f(u)}{|u|^{\alpha_{1}\alpha_{2}}}\,dx  \\
    &\leq c_{2}\int_{a}^{b}|q(x)||u|^{p-\alpha_{1}\alpha_{2}}\,dx \\
    &< c_{2}M^{p-\alpha_{1}\alpha_{2}}\int_{a}^{b}|q(x)|\,dx.
\end{align*} 
Since $\int_{a}^{b}|q(x)|\,dx < \infty$, this shows that $d$ cannot be too close to $a$ or $b$. \end{proof}

It is clear from the above proposition that we may obtain a lower bound for $||u||_{\infty}$ easily. A similar result for 
a three point condition can be formulated the same way. We present the result below. To avoid redundancy, we omit the proof.  

\begin{proposition}
Assume $(H1)-(H4)$ holds and Eq.~\eqref{1.1} has a nontrivial solution $u(x)$ satisfying \eqref{BC2}. Let $M = |u(d)| = ||u||_{\infty} = \max_{a\leq x\leq c}|u(x)| > 0$, and suppose $p-\alpha_{1}\alpha_{2} > 0.$ Then $d$ cannot be too close to $a,b,$ or $c$. 
\end{proposition} 

We now explore the changes in distance between consecutive zeros of an oscillatory solution of Eq.~\eqref{1.1}.
\begin{thm}
Assume $(H1)-(H5)$ holds and Eq.~\eqref{1.1} has a nontrivial solution $u(x)$. Suppose $\{t_{k}\}_{k=1}^{\infty}$ is an increasing sequence of zeros of $u(x)$ in $[0, \infty)$. Let there exists $\sigma > 1$ such that for any $M > 0$ 
\begin{equation}
    \int_{t}^{t+M}|q(s)|^{\sigma}\,ds \to 0 \text{ as } t \to \infty.
\end{equation} Then $t_{n+2} - t_{n} \to \infty$ as $n \to \infty$. 
\end{thm}
\begin{proof}
    Assume the contrary, that there exists $M > 0$ and a  subsequence $\{t_{n_{k}}\}_{k=1}^{\infty}$ of $\{t_{n}\}_{n=1}^{\infty}$ such that $t_{n_{k}+2} - t_{n_{k}} \leq M$ for large $k$. By the assumption, \begin{equation*}
        \int_{t_{n_{k}}}^{t_{n_{k}+2}}|q(x)|\,dx \leq \int_{t_{n_{k}}}^{t_{n_{k}}+M}|q(x)|\,dx \to 0 \text{ as } k \to \infty.
    \end{equation*} 
    By Corollary (\ref{cor2.1}), since $\psi_{1}$ is increasing we have \begin{equation}
        \int_{t_{n_{k}}}^{t_{n_{k}+2}}|q(x)|\Phi(u)\,dx > \psi_{2}\left(\frac{\frac{2}{t_{n_{k}+2} - t_{n_{k}}}}{\psi_{1}\left(\frac{t_{n_{k}+2} - t_{n_{k}}}{2}\right)}\right) \geq \psi_{2}\left(\frac{\frac{2}{M}}{\psi_{1}\left(\frac{M}{2}\right)}\right).
    \end{equation} By applying Holder's Inequality \begin{equation}
        \int_{a}^{b}|g(x)h(x)|\,dx \leq \left(\int_{a}^{b}|g(x)|^{r}\,dx\right)^{\frac{1}{r}}\left(\int_{a}^{b}|h(x)|^{s}\,dx\right)^{\frac{1}{s}},
    \end{equation} with $g(x) = q(x), h(x) = \Phi(u), r = \sigma, s = \frac{\sigma}{\sigma - 1}$, we have \begin{equation*}
        \psi_{2}\left(\frac{\frac{2}{M}}{\psi_{1}\left(\frac{M}{2}\right)}\right) < \left(\int_{a}^{b}|q(x)|^{\sigma}\,dx\right)^{\frac{1}{\sigma}}\left(\int_{a}^{b}\Phi(u)^{\frac{\sigma}{\sigma - 1}}\,dx\right)^{\frac{\sigma - 1}{\sigma}} \to 0.
    \end{equation*} This gives us a contradiction, because $\psi_{2}\left(\frac{\frac{2}{M}}{\psi_{1}\left(\frac{M}{2}\right)}\right) > 0$.
\end{proof}

In the following we discuss the zero count of a non-trivial solution of Eq.~\eqref{1.1} on a given interval.

\begin{thm}\label{thm2.3} 
Assume $(H1)-(H5)$ holds and Eq.~\eqref{1.1} has a nontrivial solution $u(x)$. 
Suppose $\{t_{k}\}_{k=1}^{2N+1}$, $N\in\mathbb{N}$ is an increasing sequence of zeroes of $u(x)$ 
in a compact interval $[a,b]$. Then 

$$
   N < \left[\psi_{2}\left(\frac{\tfrac{2}{b-a}}{\psi_{1}\left(\tfrac{b-a}{2}\right)}\right)\right]^{-1}\sum_{k=1}^{N}\max_{\xi_{k}\in [t_{2k-1},t_{2k+1}]}\left\{\int_{t_{2k-1}}^{\xi_{k}}q_{-}(x)\Phi(u)\,dx\right. \nonumber \\ 
+ \left.\int_{\xi_{k}}^{t_{2k+1}}q_{+}(x)\Phi(u)\,dx\right\},
$$
where $\Phi(u) = \frac{f(u)}{\psi_{2}(\psi_{1}(u))}$.
\end{thm}
\begin{proof}
For simplicity, we denote $L_{k}= \frac{t_{2k+1} - t_{2k-1}}{2} $.
Applying Theorem \ref{t2.2a} to the interval $[t_{2k-1}, t_{2k+1}] \subseteq [a,b]$, $k = 1, 2, \ldots, N$, we have
\begin{equation*}
    \max_{\xi_{k}\in [t_{2k-1},t_{2k+1}]}\left\{\int_{t_{2k-1}}^{\xi_{k}}q_{-}(x)\Phi(u)\,dx + \int_{\xi_{k}}^{t_{2k+1}}q_{+}(x)\Phi(u)\,dx\right\} > \psi_{2}\left(\frac{1}{L_{k}\psi_{1}\left(L_{k}\right)} \right)
\end{equation*} 
and hence
$$
    \sum_{k=1}^{N}\max_{\xi_{k}\in [t_{2k-1},t_{2k+1}]}\left\{\int_{t_{2k-1}}^{\xi_{k}}q_{-}(x)\Phi(u)\,dx + \int_{\xi_{k}}^{t_{2k+1}}q_{+}(x)\Phi(u)\,dx\right\} > \sum_{k=1}^{N}\psi_{2}\left(\frac{1}{L_{k}\psi_{1}\left(L_{k}\right)} \right).
$$
Since $\psi_{2}$ is increasing, then 
\begin{equation*}
    \psi_{2}\left(\frac{1}{L_{k}\psi_{1}(L_{k})}\right) \geq \psi_{2}\left(\frac{\tfrac{2}{b-a}}{\psi_{1}\left(\tfrac{b-a}{2}\right)}\right). 
\end{equation*} 
Therefore, we have
\begin{equation*}
    \sum_{k=1}^{N}\psi_{2}\left(\frac{1}{L_{k}\psi_{1}\left(L_{k}\right)} \right) \geq \sum_{k=1}^{N}\psi_{2}\left(\frac{\tfrac{2}{b-a}}{\psi_{1}\left(\tfrac{b-a}{2}\right)}\right) = N\psi_{2}\left(\frac{\tfrac{2}{b-a}}{\psi_{1}\left(\tfrac{b-a}{2}\right)}\right).
\end{equation*} 
Hence the proof follows.
\end{proof}

Since $q_{-}(x), q_{+}(x) \leq |q(x)|$, the following simpler result directly follows from Theorem \ref{thm2.3}.

\begin{cor}\label{cor2.2}
Let $u(x)$ be a nontrivial solution of Eq.~\eqref{1.1}. Let $\{t_{k}\}_{k=1}^{2N+1}, N \geq 1$, be an increasing sequence of zeroes of $u(x)$ in a compact interval $[a,b]$. Also, let $\psi_{2}$ be convex. Then \begin{equation*}
    N < \left[\psi_{2}\left(\frac{\frac{2}{b-a}}{\psi_{1}\left(\frac{b-a}{2}\right)}\right)\right]^{-1}\int_{a}^{b}|q(x)|\Phi(u)\,dx.
\end{equation*}
\end{cor}

\begin{rem}\label{rem2.2}
\rm{
Let $\psi_{1}$, $\psi_{2},$ and $f$ be signed power functions of the form $\psi_{1}(u) = |u|^{\alpha_{1}-1}u$, $\psi_{2}(u) = |u|^{\alpha_{2}-1}u$, and $f(u) = |u|^{\alpha_{2}\alpha_{1}-1}u$ for $\alpha_{1}, \alpha_{2} > 0$. Then Eq.~\eqref{1.1} becomes a third order quasilinear equation and in this case $\Phi(u) \equiv 1$. Therefore, by Theorem 2.4, we have 
\begin{equation*}
    \sum_{k=1}^{N}\max_{\xi_{k}\in [t_{2k-1},t_{2k+1}]}\left\{\int_{t_{2k-1}}^{\xi_{k}}q_{-}(x)\,dx + \int_{\xi_{k}}^{t_{2k+1}}q_{+}(x)\,dx\right\} > \left(\frac{2}{b-a}\right)^{\alpha_{2}(\alpha_{1} + 1)}N^{\alpha_{2}(\alpha_{1} + 1)+1}.
\end{equation*} 
This is consistent with \cite[Theorem 2.3]{R8} for third order quasilinear equations. 
Furthermore, for $\alpha_{1} = \alpha_{2} = 1$, we have by Theorem 2.4 
\begin{equation*}
    \sum_{k=1}^{N}\max_{\xi_{k}\in [t_{2k-1},t_{2k+1}]}\left\{\int_{t_{2k-1}}^{\xi_{k}}q_{-}(x)\,dx + \int_{\xi_{k}}^{t_{2k+1}}q_{+}(x)\,dx\right\} > \frac{4N^{3}}{(b-a)^{2}},
\end{equation*} which is a significant improvement on the results obtained in \cite[Theorem 6]{R7} where for arbitrarily large $N$, their result is trivial and does not provide any additional information.}
\end{rem}

\end{document}